\documentclass[12pt]{article}
\usepackage{amssymb}
\usepackage{latexsym}
\usepackage{amsmath}
\usepackage{amsthm}
\newtheorem{theorem}{Theorem}
\newtheorem{lemma}[theorem]{Lemma}
\newtheorem{proposition}[theorem]{Proposition}
\newtheorem{fact}[theorem]{Fact}

\newtheorem{conjecture}[theorem]{Conjecture}

\newtheorem{question}[theorem]{Question}

\begin{document}
\title{On a Question of Abraham Robinson's}
\author{Jochen Koenigsmann\\ Accepted for publication in Israel Journal of Mathematics}
\date{}
\maketitle
\begin{abstract}
In this note we give a negative answer to Abraham Robinson's question
whether a finitely generated extension of an undecidable field
is always undecidable.
We construct `natural' undecidable fields of transcendence degree $1$ over $\mathbb{Q}$ all of whose proper finite extensions are decidable.
We also construct undecidable algebraic extensions of $\mathbb{Q}$
that allow decidable finite extensions.\footnote{\noindent The research
on this note was undertaken
while the author enjoyed the hospitality of
the American Institute of Mathematics (AIM)
which has been greatly appreciated.
The author would also like to thank Carlos Videla for pointing him
to Robinson's question.}
\end{abstract}
A field $F$ is said to be {\em decidable}
if {\tt Th}$(F)$, the first order theory of $F$
in the language ${\cal L}:=\{+,\times;0,1\}$ of rings, is decidable.
Equivalently, $F$ is decidable
if there is an effective axiomatization of {\tt Th}$(F)$,
i.e., an algorithm producing a (typically infinite) list of 
${\cal L}$-sentences true in $F$ from which every
sentence in {\tt Th}$(F)$ can be deduced.

It is straightforward to see that
if $F$ is decidable and $E$ is a finite extension of $F$
obtained by adjoining elements algebraic over the prime field
that then $E$ is decidable as well.

40 years ago,
in his retiring presidential address [R] at the annual ASL meeting in Dallas,
Abraham Robinson asked the following question in the opposite direction:
{\em Is a finitely generated extension of an undecidable field always undecidable?}
It turns out that the answer to this question is `No' in a very strong sense: 
\begin{theorem}{\em (A variation of Example 4.4 in [CDM]\footnote{A week after presenting this Theorem at AIM we learned (via private communication)
that a very similar result had already been proved in the unpublished manuscript [CDM],
where the existence of an `outlandish example' of an undecidable field $F$ of infinite transcendence degree over $\mathbb{Q}$ with all proper finite extensions being decidable and isomorphic has been established, though in a less constructive way.} )}
There are uncountably many pairwise non-elementarily-equivalent undecidable algebraic extensions $F$ of $\overline{\mathbb{Q}}(t)$
for which every proper finite extension $E/F$ is decidable
(and all these $E$ are elementarily equivalent).
\end{theorem}
Here $\overline{\mathbb{Q}}$ denotes the algebraic closure of $\mathbb{Q}$
and $\overline{\mathbb{Q}}(t)$ the rational function field in one variable
over $\overline{\mathbb{Q}}$.

In the light of this theorem it is natural to narrow down
Robinson's question to fields that are algebraic over $\mathbb{Q}$,
that is, to ask
whether every finite extension of an undecidable algebraic extension $F$ of $\mathbb{Q}$ is always undecidable
(a question also raised in the concluding Remark of [D]).
This is known to be the case if $F$ is a number field,
i.e., a finite extension of $\mathbb{Q}$:
all number fields are undecidable.
In general, however, the answer is again in the negative:
\begin{theorem}
For every prime $p$ there are infinitely many pairwise non-isomorphic
(and hence, here, non-elementarily-equivalent) decidable algebraic extensions
$E$ of $\mathbb{Q}$
having uncountably many undecidable subfields $F$ of degree $[E:F] = p$.
\end{theorem}

While Robinson had asked his question about {\em finitely generated}
field extensions $E/F$
we have, in both of the above theorems, even produced examples of
{\em finite}, i.e., finitely generated {\em algebraic} field extensions $E/F$
such that $E$ was decidable and $F$ undecidable.
As it happens, these examples were bound to be of this kind,
as Robinson's question when asked about finitely generated {\em non-algebraic} field extensions is very likely to have a positive answer
in the light of the following well supported:
\begin{conjecture}
If $F$ is an arbitrary field and $E/F$ a finitely generated non-algebraic field extension then $E$ is undecidable.
\end{conjecture}
Note that the conjecture does not assume $F$ to be undecidable.

This conjecture which has been raised by Malcev in [M] in the case that $E$ is a rational function field (so purely transcendental over $F$) has, for example, been proved in the following cases:
\begin{itemize}
\item
if the characteristic of $F$ is $>0$ ([ES])
\item
if $E$ is the function field of a curve over $F$,
if $F$ is undecidable and either $F$ is large
(in the technical sense that varieties over $F$ with one $F$-rational point
have infinitely many) or, for some $n>1$, $F^\times /(F^\times)^n$ is finite:
in this case, $F$ is definable in $E$ (Theorem 2 of [K])
\item
if $E=F(t)$ and $F$ is formally real ([RR])
\item
if $E=F(t)$, $F$ is infinite and existentially undecidable
(for then {\tt Th}$_\exists(E)=$ {\tt Th}$_\exists(F)$)
\end{itemize}
The most irritating unknown case is $\mathbb{C}(t)$: it is not known
whether or not $\mathbb{C}(t)$ is decidable.

As mentioned in the beginning,
if $E/F$ is a finite extension that can be
generated by elements algebraic over the prime field
then decidability goes up (from $F$ to $E$).
For arbitrary finite extensions $E/F$, however,
in order to interpret $E$ in $F^n$ (where $n=[E:F]$)
one needs parameters for the minimal polynomials of
a basis for $E$ over $F$,
yet it can happen that a field $F$ is decidable,
but becomes undecicable if one allows new constants:
For example if $F=\mathbb{R}$, any transcendental
element $c\in\mathbb{R}$ which realizes an undecidable Dedekind cut, 
makes the theory with a constant for $c$ undecidable.
We have no answer to the following twist of Robinson's question:

\begin{question}
Is a finite extension of a decidable field always decidable?
\end{question}

When Robinson presented his `metamathematical problems'
from twelve different areas of mathematics
(the problem addressed in this note is only one of them)
he emphasized that these problems seem of interest
`{\em not only for their own sake
but also because their solution might well require weapons
whose introduction would close definite gaps in our armory}'.
The new weapons used in this note for solving the problem
are, indeed, several deep results from field arithmetic,
an area of mathematics that had rapidly developed only after Robinson's death in 1974.

For the convenience of the reader, in section 1,
we shall recall these results
which have by now become classics in their own right.
Sections 2 and 3 are then devoted to proving
Theorem 1 and 2 respectively.
\medskip\\
{\bf Acknowledgement:}
The author would like to thank Moshe Jarden
for suggesting very helpful improvements on an earlier draft.
In particular, the proof of Proposition \ref{PACdec} is due to him.
%%%%%%%%%%%%%%%%%%%%%%%%%%%%%%%%%%%%%%%%%%%%%%%%%%%%%%%%%%%%%%%%
\section{A few classics from field arithmetic}
Let us recall that a field $K$ is called {\bf Hilbertian}
if it satisfies Hilbert's Irreducibility Theorem,
i.e., if for every irreducible polynomial $f(T,X)\in K[T,X]$
there are infinitely many $t\in K$ such that $f(t,X)\in K[X]$ is irreducible
(this property goes up to finite extensions of $K$).

$K$ is called {\bf PAC} (pseudo algebraically closed)
[{\em resp. {\bf PRC} (pseudo real closed)}]
if every non-empty variety $V$ defined over $K$
has a $K$-rational point
[{\em provided it has rational points in all real closures of $K$}].
Both properties are inherited by any algebraic extension of $K$.
If $K$ is not formally real then PRC and PAC are the same.

We denote the separable algebraic closure of $K$ by $\overline{K}$
and the absolute Galois group {\tt Gal}$(\overline{K}/K)$ by $G_K$.

Finally, we denote by $\hat{F}_\omega$ the free profinite group
on countably infinitely many generators.
\begin{fact}{\em ([FV], also Example 24.8.5(b) in [FJ])}\\
\label{FV}
If $K$ is a countable Hilbertian PAC field then $G_K\cong \hat{F}_\omega$.
\end{fact}
Conversely, it is also true that if $K$ is a countable PAC field with 
$G_K\cong \hat{F}_\omega$ then $K$ is Hilbertian
(a result of Roquette, cf. Corollary 27.3.3 in [FJ])),
a fact that we shall not use.
What we will use, however, is the following celebrated result
of Roquette's student Weissauer:
\begin{fact}{\em ([W], Satz 9.7, also Theorem 13.9.1(b) in [FJ])}\\
\label{Weissauer}
Let $K$ be a Hilbertian field, let $F$ be a Galois extension of $K$
and let $F^\prime$ be a proper finite separable extension of $F$.
Then $F^\prime$ is Hilbertian.
\end{fact}
\begin{fact}{\em (A special case of a geometric local-global principle
proved by Moret-Bailly in [MB], with a more elementary proof in [GPR])}\\
\label{GPR}
Let $K$ be a countable Hilbertian field, let ${\cal P}$ be a non-empty finite set of orderings of $K$ and assume that
either $K$ admits an ordering not in ${\cal P}$
or, for some prime $p$, $K$ can be embedded into $\mathbb{Q}_p$, the field of $p$-adic numbers.
Let $K^{\cal P}$ be the intersection of all real closures of $K$ w.r.t. all orderings in ${\cal P}$.
Then $K^{\cal P}$ is PRC and hence, by Fact \ref{Weissauer},
$K^{\cal P} (\sqrt{-1})$ is Hilbertian and PAC. 
\end{fact}
For the following characterization of $\hat{F}_\omega$ due to Iwasawa
define a {\bf finite embedding problem} for the profinite group $G$
to be a pair of (continuous) epimorphisms $(\phi:\; G\to A,\;\alpha:\; B\to A)$
where $B$ (and hence $A$) is a finite group.
A {\bf solution} of the embedding problem is an epimorphism $\gamma:\; G\to B$
such that $\alpha\circ\gamma = \phi$
(sometimes these are called {\em proper} solutions
while `solutions' are not assumed to be onto).
\begin{fact}{\em ([I], p. 567, also generalized in Theorem 24.8.1 in [FJ])}\\
\label{Iwasawa}
Let $G$ be a profinite group of at most countable rank.
Then $G\cong \hat{F}_\omega$
if and only if every finite embedding problem for $G$ is solvable.
\end{fact}  
For a field $F$ of characteristic $0$ we let
$F^{alg}:=F\cap\overline{\mathbb{Q}}$ denote its algebraic part
and we let {\tt Th}$^{alg}(F)$ be the subset of {\tt Th}$(F)$
that determines $F^{alg}$ (up to isomorphism) - it is axiomatized
by saying which (irreducible, monic) polynomials in $\mathbb{Z}[X]$ have a zero in $F$ and which ones don't.
Note that {\tt Th}$^{alg}(F)$,
while determining $F^{alg}$ up to isomorphism,
is still not an axiomatization of {\tt Th}$(F^{alg})$:
{\tt Th}$^{alg}(F)$ is never complete.

\begin{proposition}
\label{PACdec}
Let $F$ be a countable PAC field of characteristic $0$
with $G_F\cong\hat{F}_\omega$.
Then $F$ is decidable if and only if {\tt Th}$^{alg}(F)$ is decidable.
\end{proposition}
{\em Proof:}
We extend the language ${\cal L}=\{ +,\times;0,1\}$ of rings
by predicates $R_n$, one for each $n\in\mathbb{N}$.
Let $T_{R,0}$ be the theory of fields of characteristic $0$ 
with algebraic part $F^{alg}$ that satisfy all of the axioms
$$ R_n(X_0,\ldots ,X_{n-1})\longleftrightarrow
(\exists Z)[Z^n + X_{n-1}Z^{n-1} +\cdots + X_0 = 0].$$
Let $L/K$ be an extension of fields of characteristic $0$
with algebraic part $F^{alg}$.
Then $L/K$ is also an extension of models of $T_{R,0}$
if and only if $K$ is algebraically closed in $L$.
Since we are in characteristic $0$,
the latter condition means that $L$ is a regular extension of $K$.

By a small variation of Theorem 27.2.3 of [FJ],
$T_{R,0}$ has a model companion $\widetilde{T}_{R,0}$
whose models are $\omega$-free PAC-fields of characteristic $0$
with algebraic part $F^{alg}$
(`$\omega$-free' means that any countable elementary substructure
has absolute Galois group $\hat{F}_\omega$).
By the preceding paragraph,
$\widetilde{T}_{R,0}$ has the amalgamation property.
Hence $\widetilde{T}_{R,0}$ is even a completion of $T_{R,0}$.
By Satz 3.22 of [Pr], it follows that $\widetilde{T}_{R,0}$ admits elimination of quantifiers.
Therefore, in order to decide
whether a sentence in our extended language holds in $F$,
we have to be able to decide which sentences of the form
$R_n(a_0,\ldots ,a_{n-1})$ with $a_0,\ldots ,a_{n-1}\in\mathbb{Z}$
hold in $F$,
i.e., which polynomials $Z^n+a_{n-1}Z^{n-1}+\cdots +a_0$ have a root in $F^{alg}$.
This is doable if and only if 
{\tt Th}$^{alg}(F)$ is decidable.
\qed
%%%%%%%%%%%%%%%%%%%%%%%%%%%%%%%%%%%%%%%%%%%%%%%%%%%%%%%%%%%%%%%%
\section{Proof of Theorem 1}
Let $K$ be an algebraic extension of $\overline{\mathbb{Q}}(t)$
which is Hilbertian and PAC.
Such extensions exist:
for example, let $R:=\mathbb{R}\cap\overline{\mathbb{Q}}$
be the field of real algebraic numbers and let
${\cal P}$ be a non-empty finite set of orderings on $R(t)$.
Then, by Fact \ref{GPR}, the field $K:=R(t)^{\cal P}(\sqrt{-1})$
is Hilbertian, PAC and, by construction, an algebraic extension
of $\overline{\mathbb{Q}}(t)$
(note that $R(t)$ allows infinitely many orderings).\footnote{Alternatively,
one could have used the fact that,
(in the sense of the Haar measure on profinite groups)
for almost all $\sigma\in G_{\overline{\mathbb{Q}}(t)}$,
the maximal Galois extension $K$ of $\overline{\mathbb{Q}}(t)$
with $\sigma\in G_K$ is PAC and Hilbertian
(Theorem 2.7 in [J], also Theorem 27.4.8 in [FJ]),
but this is certainly less constructive.}

Let $\Sigma$ be any set of primes and let $K(\Sigma)$
be the compositum of all $K(p)$ for $p\in\Sigma$
(inside a fixed algebraic closure of $K$),
where $K(p)$ denotes the maximal pro-$p$ Galois extension of $K$.

Now construct a chain of algebraic extensions of $K$
$$K_0\subseteq K_1\subseteq K_2\subseteq\ldots,$$
where $K_0:=K$ and, for $n>0$, $K_n:= K_{n-1}(\Sigma )$.
Then $F:=K_\infty (\Sigma):= \bigcup_{n=0}^\infty K_n$
is a Galois extension of $K$:
note that $G_{K_n}$ is a characteristic subgroup of $G_{K_{n-1}}$,
so all $G_{K_n}$ are normal (and, in fact, characteristic) subgroups of $G_K$.
And, for any prime $p$,
$$F=F(p)\Longleftrightarrow p\in\Sigma:$$
`$\Leftarrow$' follows by construction.
For `$\Rightarrow$', assume $p\not\in\Sigma$.
As $K$ is Hilbertian, $K\neq K(p)$ and, by construction,
$K(p)\cap F=K$. So $F\neq FK(p)\subseteq F(p)$.

In particular, if $\Sigma\neq\Sigma^\prime$
then $K_\infty (\Sigma)$ and $K_\infty (\Sigma^\prime)$
are not elementarily equivalent, so there are
uncountably many pairwise non-elementarily-equivalent such $F$,
of which only countably many can be decidable.

By Weissauer's Theorem (Fact \ref{Weissauer}),
every proper finite extension $E/F$ is again Hilbertian
(and as algebraic extension of the PAC field $K$ still PAC).
Hence, by Fact \ref{FV}, $G_E\cong \hat{F}_\omega$.
Moreover, by construction, $E^{alg}=\overline{\mathbb{Q}}$,
so {\tt Th}$^{alg}(E)$ is decidable
(all $f\in\mathbb{Z}[X]\setminus\mathbb{Z}$ have zeros)
and thus, by Proposition \ref{PACdec}, so is $E$.
\qed

If one doesn't like the counting argument
one may, alternatively, choose a set $\Sigma$ of primes
which is not recursive in order to guarantee that
$K_\infty (\Sigma)$ is not decidable.
%%%%%%%%%%%%%%%%%%%%%%%%%%%%%%%%%%%%%%%%%%%%%%%%%%%%%%%%%%%%%%%%
\section{Proof of Theorem 2}
Let us first prove the following 
\begin{lemma}
\label{cod1}
Let $V$ be a vector space of infinite dimension over a field $k$,
let $G$ be a (profinite) group acting on $V$ as group of $k$-linear transformations,
asssume that all orbits of $V$ under this action of $G$ are finite
and that, for every finite-dimensional $G$-invariant subspace $W$ of $V$,
there is a non-zero $G$-invariant subspace $W^\prime$ of $V$
with $W\cap W^\prime = \{ 0\}$.
Then, under the induced action of $G$ on the set $\widetilde{V}$
of all codimension $1$ subspaces of $V$,
$\widetilde{V}$ has uncountably many $G$-orbits.
\end{lemma}
{\em Proof:}
Choose any $b_1\in V\setminus \{ 0\}$ and let
$B_1$ be a maximal linearly independent subset of the (finite) $G$-orbit of $b_1$, so every element in the $G$-orbit of $b_1$ is in $\langle B_1\rangle$,
the span of $B_1$.
Now assume we have constructed finite sets $B_1,B_2,\ldots, B_m$
such that each $B_i$ spans the whole $G$-orbit of any of its members,
$\bigcup_{i=1}^m B_i$ is linearly independent,
and $W:=\langle B_1\cup\ldots \cup B_m\rangle = \langle B_1\rangle\oplus\cdots\oplus\langle B_m\rangle$.
Then, by assumption, there is a $G$-invariant non-zero subspace $W^\prime$ with $W\cap W^\prime =\{ 0\}$.
So we can choose $b_{m+1}\in W^\prime\setminus\{ 0\}$
and let $B_{m+1}$ be a maximal linearly independent subset of the $G$-orbit of $b_{m+1}$.
Then $B_1\cup\ldots\cup B_{m+1}$ is linearly independent and
$\langle B_1\cup\ldots\cup B_{m+1}\rangle = B_1\oplus\cdots\oplus B_{m+1}$.

Now prolong the linearly independent set $\bigcup_{i=1}^\infty B_i$ to a basis $B$ of $V$ and define, for each proper subset $I$ of $\mathbb{N}$,
a map $\chi_I:B\to k$ by setting
$$\chi_I(b)=\left\{\begin{array}{ll}
0 & \mbox{if }b\in\bigcup_{i\in I}B_i\\
1 & \mbox{if }b\in B\setminus\bigcup_{i\in I}B_i
\end{array}\right.$$
Let $T_I:V\to k$ be the $k$-linear map extending $\chi_I$.
As $I$ was a proper subset of $\mathbb{N}$, the image of $T_I$ is $k$,
and so $\ker T_I$ is a codimension 1 subspace of $V$,
so $\ker T_I\in\widetilde{V}$.
Moreover, for $I\neq I^\prime$, 
$\ker T_I\neq\ker T_{I^\prime}$ and, stronger still,
$\ker T_I$ and $\ker T_{I^\prime}$ are in different $G$-orbits:
otherwise there is some $\sigma\in G$ such that
$$\ker T_{I^\prime} = \{\sigma (v)\mid v\in\ker T_I\}.$$
But then, say for $i\in I\setminus I^\prime$ and for any $b\in B_i$,
$b\in\ker T_I$, and, since, by construction, $\sigma (b)\in\langle B_i\rangle$,
also $\sigma (b)\in\ker T_I$, whereas $T_{I^\prime}(b)=1$: contradiction.

As there are uncountably many proper subsets $I$ of $\mathbb{N}$,
this shows that $\widetilde{V}$ has uncountably many distinct $G$-orbits.
\qed 
\medskip

To prove Theorem 2,
let $p$ be a given fixed prime.
Let $K=\mathbb{Q}$ and let ${\cal P}$
be the singleton set containing the unique ordering on $\mathbb{Q}$.
As $\mathbb{Q}\subseteq\mathbb{Q}_l$
for some (in fact, all) primes $l$,
$(K,{\cal P})$ satisfies the hypotheses of Fact \ref{GPR}
and hence the field $L:=K^{\cal P}(\sqrt{-1})$ is Hilbertian and PAC
(the field $K^{\cal P}=\mathbb{Q}^{tot.r.}$
is called the field of totally real numbers,
it is the maximal Galois extension of $\mathbb{Q}$ inside $\mathbb{R}$).
Now let
$$E = L(\sqrt[p]{L}):=L(\{ \sqrt[p]{a}\mid a\in L\}).$$
Then $E/L$ is a Galois extension
($L$ contains a primitive $p$-th root $\zeta_p$ of unity,
as $\zeta_p +\zeta_p^{-1}\in\mathbb{Q}^{tot.r.}$
and $\sqrt{4-(\zeta_p+\zeta_p^{-1})^2}\in\mathbb{Q}^{tot.r.}$).
As $L$ is Hilbertian, $L^\times /(L^\times)^p$ is an
infinite-dimensional $\mathbb{F}_p$-vector space,
and, for any codimension 1 subspace $U$ of $L^\times /(L^\times)^p$,
the field $F_U:=L(\sqrt[p]{U})$ is a Galois subextension of $E/L$
with $[E:F]=p$.
Hence, by Fact \ref{Weissauer}, $E$ is Hilbertian, and still PAC.
By Fact \ref{FV}, this implies that $G_E\cong \hat{F}_\omega$.
Further note that {\tt Th}$^{alg}(E)$ is decidable:
For each $n>0$ and each $\overline{a}=(a_0,\ldots ,a_{n-1})\in\mathbb{Z}^n$
let $$f_{\overline{a}}(X):=X^n + a_{n-1}X^{n-1}+\ldots + a_1X + a_0.$$
Then the sets
\begin{itemize}
\item $I_n:=\{ \overline{a}\in\mathbb{Z}^n\mid f_{\overline{a}}(X)\mbox{ is irreducible}\}$,
\item $T_n:=\{ \overline{a}\in I_n\mid f_{\overline{a}}\mbox{ splits in }\mathbb{Q}^{tot.r}\}$
(that the splitting field of $f_{\overline{a}}$ is totally real is equivalent to $f_{\overline{a}}$ having $n$ distinct roots in $\mathbb{R}$
which, by (effective!) quantifier elimination for $\mathbb{R}$ in ${\cal L}\cup\{ <\}$,
is equivalent to some polynomial equations and inequalities in the coefficients $\overline{a}$), and
\item ${\cal E}_n$ of those $\overline{a}\in I_n$
for which the splitting field of $f_{\overline{a}}$ is contained
in an elementary abelian $p$-Galois extension of $M(\sqrt{-1},\zeta_p)$
for some totally real Galois extension $M/\mathbb{Q}$
with $[M:\mathbb{Q}]\leq n!$
\end{itemize}
are all decidable (= recursive) and hence, so is {\tt Th}$^{alg}(E)$,
being axiomatized by
$$\begin{array}{rl}
 & \{\neg\exists x\; f_{\overline{a}}(x)=0\vee\exists x_1,\ldots ,x_n\;
f_{\overline{a}}(X)=\prod_{i=1}^n (X-x_i)\mid \overline{a}\in I_n, n\in\mathbb{N}\}\\
\cup & \{\neg\exists x\; f_{\overline{a}}(x)=0\mid \overline{a}\in I_n\setminus{\cal E}_n, n\in\mathbb{N}\}\\
\cup & \{\exists x_1,\ldots ,x_n\; f_{\overline{a}}(X)=\prod_{i=1}^n(X-x_i)\mid \overline{a}\in {\cal E}_n, n\in\mathbb{N}\}
\end{array}.$$
Hence, by Proposition \ref{PACdec}, $E$ is decidable.

To see that there are uncountably many $U$ as above
for which the corresponding $F_U$ are pairwise non-elementarily equivalent,
let $G:=\mbox{\tt Gal}(L/\mathbb{Q})$
and consider the natural action of $G$
on the $\mathbb{F}_p$-vectorspace $V:=L^\times/(L^\times)^p$:
$$\begin{array}{ccc}
G\times L^\times /(L^\times )^p & \longrightarrow & L^\times /(L^\times )^p\\
(\sigma , a (L^\times)^p) & \longmapsto & (\sigma a)(L^\times )^p
\end{array}$$
As each $a\in L^\times$ lies in a finite Galois subextension $L_a/\mathbb{Q}$ of $L/\mathbb{Q}$,
each $G$-orbit of $L^\times /(L^\times)^p$ is finite:
the orbit of $a(L^\times)^p$ is $\{\tau (a)\cdot (L^\times )^p\mid\tau\in\mbox{\tt Gal}(L_a/\mathbb{Q})\}$.

Moreover, for any finite(-dimensional) $G$-invariant subspace $W$ of $V$,
there is a non-zero $G$-invariant subspace $W^\prime$ of $V$
with $W\cap W^\prime =\{ 0\}$:
If $p>2$ then the canonical $\mathbb{F}_p$-vector space homomorphisms
$$\mathbb{Q}^\times /(\mathbb{Q}^\times)^p \longrightarrow \mathbb{Q}^{tot.r.\times} / (\mathbb{Q}^{tot.r.\times})^p \longrightarrow
L^\times /(L^\times)^p$$
are all injective,
so, for any $a\in\mathbb{Q}^\times\setminus (\mathbb{Q}^\times )^p$ with $a(L^\times)^p\not\in W$,
$W^\prime:=\langle a(L^\times)^p\rangle$ is a $1$-dimensional $G$-invariant subspace of $V$ with $W\cap W^\prime =\{ 0\}$.
If $p=2$, then for any prime $q$, $\sqrt{q}$ spans a $1$-dimensional $G$-invariant subspace of $V$:
$\sqrt{q}(L^\times)^2 = (-\sqrt{q})(L^\times)^2$,
so, again, we find $W^\prime$ as required:
note that $\sqrt{q}\not\in\pm(\mathbb{Q}^{tot.r.\times})^2$
and that the kernel of
$\mathbb{Q}^{tot.r.\times} / (\mathbb{Q}^{tot.r.\times})^2 \longrightarrow
L^\times /(L^\times)^2$ is
$\pm(\mathbb{Q}^{tot.r.\times})^2 / (\mathbb{Q}^{tot.r.\times})^2$.

Now, by Lemma \ref{cod1},
there are uncountably many codimension 1 subspaces $U$ of $V$
all lying in distinct $G$-orbits.
But for any two codimension 1 subspaces $U,U^\prime$ of $V$ one has
$$\exists\sigma\in G\;U^\prime =\sigma(U)\Longleftrightarrow F_{U^\prime}\cong F_U\Longleftrightarrow F_{U^\prime}\equiv F_U,$$
so that there are uncountably many pairwise non-elementarily equivalent such $F_U$, as claimed.

This gives just one $E$ of the kind Theorem 2 wants to have it.
To get infinitely many, let $l$ be any prime $\neq 2,p$
and let $E_l:=E(\sqrt[l]{2})$.
Then $E_l$ (being parameter-free interpretable in $E^l$) is still decidable
and there are uncountably many undecidable $F_{U,l}:=F_U(\sqrt[l]{2})$
with $U$ as above and $[E_l:F_{U,l}]=p$.
Finally, observe that $\sqrt[l]{2}\not\in\mathbb{Q}^{tot.r.}$ (as $l>2$),
so $\sqrt[l]{2}\not\in E$ (as $l\neq p$) and, for $l\neq l^\prime$,
$E_l\not\equiv E_{l^\prime}$.\qed
%%%%%%%%%%%%%%%%%%%%%%%%%%%%%%%%%%%%%%%%%%%%%%%%%%%%%%%%%%%%%%%%

Mathematical Institute, Radcliffe Observatory Quarter, Woodstock Road, Oxford OX2 6GG, UK\\
{\tt koenigsmann@maths.ox.ac.uk}
\end{document}